\documentclass[reqno,12pt]{amsart}

\usepackage{epsf}
\usepackage{graphics}
\usepackage{amssymb}
\usepackage{amsmath}

\date{}

\theoremstyle{plain}
\newtheorem{theorem}{Theorem}
\newtheorem{corollary}{Corollary}

\theoremstyle{definition}

\theoremstyle{remark}

\def\N{{\mathbb N}}

\def\R{{\mathbb R}}

\title{Positive braids of maximal signature} 

\author{Sebastian Baader}

\begin{document}

\begin{abstract} We characterise positive braid links with positive Seifert form via a finite number of forbidden minors. From this we deduce a one-to-one correspondence between prime positive braid links with positive Seifert form and simply laced Dynkin diagrams, as well as a simple classification of alternating positive braid knots. 
\end{abstract}

\maketitle

\section{Introduction}

A minor of a graph is obtained by deleting and contracting a finite number of edges of that graph. The Wagner conjecture states that a family of mutually non-minor graphs is necessarily finite~\cite{L}. In particular, every minor-closed property of graphs is characterised by a finite number of forbidden minors. Thinking of graphs as spines of surfaces embedded in $\R^3$, we are naturally led to define a minor of a Seifert surface as an incompressible subsurface (i.e. a subsurface whose complement has no disc component). Here a Seifert surface is a  compact orientable surface with boundary smoothly embedded in $\R^3$, considered up to isotopy. We observe right away that infinite families of mutually non-minor Seifert surfaces exist, e.g. families of embedded annuli with pairwise non-isotopic core curves.

In this paper we study a natural minor-closed property derived from the Seifert form. We say that the Seifert form of a Seifert surface is positive, if its symmetrisation ($V+V^{\mathrm{T}}$ in matrix notation) is positive-definite. In the case of fibred links, there is a canonical Seifert form associated with the fibre surface. This applies to positive braid links.

\begin{theorem} Let $L \subset \R^3$ be a positive braid link with non-positive Seifert form. Then the fibre surface of $L$ contains an incompressible subsurface of type $T$, $E$, $X$ or $Y$.
\end{theorem}

Here the surfaces of type $T$, $E$, $X$ and $Y$ are the canonical fibre surfaces associated with the positive braids $\sigma_1^4 \sigma_2 \sigma_1^4 \sigma_2$, $\sigma_1^6 \sigma_2 \sigma_1^3 \sigma_2$, $\sigma_1 \sigma_2^2 \sigma_1 \sigma_3 \sigma_2^2 \sigma_3$ and $\sigma_1^3 \sigma_2^2 \sigma_1^3 \sigma_2$, respectively. The meaning of these will shortly become clear. The above criterion has a remarkable consequence.

\begin{theorem} There is a natural one-to-one correspondence between prime positive braid links with positive Seifert form and simply laced Dynkin diagrams, reading as follows:

\smallskip
$A_n$ \quad $\sigma_1^n$

\smallskip
$D_n$ \quad $\sigma_1^n \sigma_2 \sigma_1^2 \sigma_2$

\smallskip
$E_6$ \quad $\sigma_1^3 \sigma_2 \sigma_1^3 \sigma_2$

\smallskip
$E_7$ \quad $\sigma_1^4 \sigma_2 \sigma_1^3 \sigma_2$

\smallskip
$E_8$ \quad $\sigma_1^5 \sigma_2 \sigma_1^3 \sigma_2$
\end{theorem}

Naturality means that the quadratic form associated with the Dynkin diagrams is isomorphic to the symmetrised Seifert form of the corresponding braids. The knots arising from the series $A_n$ and $E_6$, $E_8$ are the torus knots of type $T(2,n)$ ($n \in \N$ odd), $T(3,4)$ and $T(3,5)$, respectively. The other diagrams give rise to links with two or three components. 

Positive braid links are known to have positive signature invariant~\cite{Ru}. Positivity of the Seifert form is a much stronger feature. For positive braid knots it amounts to the equality $\sigma=2g$ between the signature invariant and the minimal genus of the knot. There is an analogous equality $s=2g$ between Rasmussen's invariant and the minimal genus, which holds for all positive braid knots~\cite{Ra}. We conclude that positive braid knots with non-positive Seifert form are non-alternating (since alternating knots satisfy $s=\sigma$). We are left with $T(2,n)$, $T(3,4)$ and $T(3,5)$ as potential prime alternating positive braid knots. The latter two are non-alternating.

\begin{corollary} The torus knots of type $T(2,n)$ are the only prime alternating positive braid knots.
\end{corollary}  

According to Cromwell~\cite{C}, the prime components of a positive braid link are positive braid links. Therefore alternating positive braid knots are connected sums of torus knots of type $T(2,n)$. This statement was previously proved by Nakamura within the framework of homogeneous links~\cite{N}.


\section{Positive braid links}

A positive braid link is the closure of a positive braid, i.e. a finite product of the standard braid group generators $\sigma_1, \sigma_2 \ldots$. It is well-known that positive braid links are fibred; their fibre surface is the standard Seifert surface associated with the braid diagram~\cite{St}. We use brick diagrams as a suggestive notation for positive braids and their fibre surfaces. Brick diagrams for the braids appearing in Theorems~1 and~2 are depicted in Figure~1. The small rectangles of a brick diagram (with the positive orientation) can be thought of as a basis for the first homology group of the fibre surface, since the latter retracts on the brick diagram.

\begin{figure}[ht]
\scalebox{1.2}{\raisebox{-0pt}{$\vcenter{\hbox{\epsffile{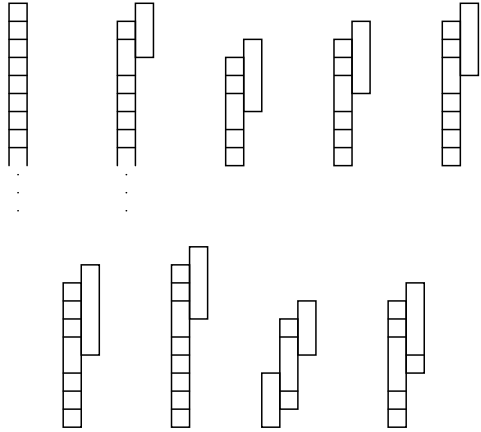}}}$}} 
\caption{Brick diagrams for $A_n$, $D_n$, $E_6$, $E_7$, $E_8$ and $T$, $E$, $X$, $Y$}
\end{figure}

It is easy to read off the Seifert form from a brick diagram. In the above examples the rectangles are linked in a tree-like pattern (see Figure~2). Here two rectangles are linked if and only if they are arranged as in the braids $\sigma_1^3$, $\sigma_1 \sigma_2 \sigma_1 \sigma_2$ or $\sigma_2 \sigma_1 \sigma_2 \sigma_1$. Indeed, these all represent the positive trefoil knot. The rectangles of the braid $\sigma_1 \sigma_2 \sigma_2 \sigma_1$ are not linked since this represents a connected sum of two positive Hopf links. 

\begin{figure}[ht]
\scalebox{1.2}{\raisebox{-0pt}{$\vcenter{\hbox{\epsffile{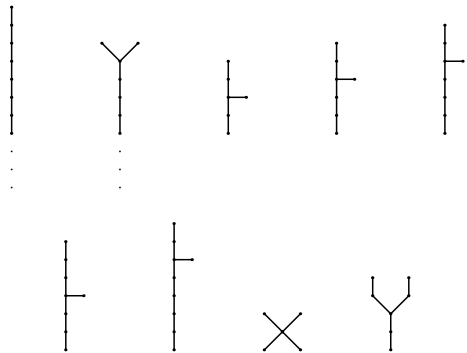}}}$}} 
\caption{Linking trees for $A_n$, $D_n$, $E_6$, $E_7$, $E_8$ and $T$, $E$, $X$, $Y$}
\end{figure}

The resulting symmetrised Seifert matrices are nothing but the Cartan matrices of the linking trees, viewed as Coxeter systems. The first five are positive definite since they correspond to spherical Coxeter systems. The remaining four matrices are
$$
\begin{pmatrix}
2 & 1 & 1 & 0 & 0 & 1 & 0 & 0 \\
1 & 2 & 0 & 0 & 0 & 0 & 0 & 0 \\
1 & 0 & 2 & 1 & 0 & 0 & 0 & 0 \\
0 & 0 & 1 & 2 & 1 & 0 & 0 & 0 \\
0 & 0 & 0 & 1 & 2 & 0 & 0 & 0 \\
1 & 0 & 0 & 0 & 0 & 2 & 1 & 0 \\
0 & 0 & 0 & 0 & 0 & 1 & 2 & 1 \\
0 & 0 & 0 & 0 & 0 & 0 & 1 & 2 \\
\end{pmatrix},
\begin{pmatrix}
2 & 1 & 1 & 0 & 1 & 0 & 0 & 0 & 0 \\
1 & 2 & 0 & 0 & 0 & 0 & 0 & 0 & 0 \\
1 & 0 & 2 & 1 & 0 & 0 & 0 & 0 & 0 \\
0 & 0 & 1 & 2 & 0 & 0 & 0 & 0 & 0 \\
1 & 0 & 0 & 0 & 2 & 1 & 0 & 0 & 0 \\
0 & 0 & 0 & 0 & 1 & 2 & 1 & 0 & 0 \\
0 & 0 & 0 & 0 & 0 & 1 & 2 & 1 & 0 \\
0 & 0 & 0 & 0 & 0 & 0 & 1 & 2 & 1 \\
0 & 0 & 0 & 0 & 0 & 0 & 0 & 1 & 2 \\
\end{pmatrix},
$$
$$
\begin{pmatrix}
2 & 1 & 1 & 1 & 1 \\
1 & 2 & 0 & 0 & 0 \\
1 & 0 & 2 & 0 & 0 \\
1 & 0 & 0 & 2 & 0 \\
1 & 0 & 0 & 0 & 2 \\
\end{pmatrix},
\begin{pmatrix}
2 & 1 & 0 & 1 & 0 & 1 & 0  \\
1 & 2 & 1 & 0 & 0 & 0 & 0  \\
0 & 1 & 2 & 0 & 0 & 0 & 0  \\
1 & 0 & 0 & 2 & 1 & 0 & 0  \\
0 & 0 & 0 & 1 & 2 & 0 & 0  \\
1 & 0 & 0 & 0 & 0 & 2 & 1  \\
0 & 0 & 0 & 0 & 0 & 1 & 2  \\
\end{pmatrix}
$$
and have determinant zero
\footnote{Incidentally, these correspond to the simply laced affine Coxeter systems (see~\cite{AC} for an elegant characterisation of spherical and affine Coxeter systems).}.
As a consequence, every Seifert surface that contains one of the surfaces $T$, $E$, $X$, $Y$ as a minor has non-positive Seifert form. We will prove the converse in the next section.

\section{Classification}

The proofs of Theorems~1 and~2 rely on careful considerations on $3$- and $4$-stranded positive braids.

\medskip
\noindent
\textbf{$3$-braids}

\noindent
Let $L$ be a link represented by a positive $3$-braid $\beta$. Applying the braid relation $\sigma_2 \sigma_1 \sigma_2=\sigma_1 \sigma_2 \sigma_1$ and conjugation to $\beta$ does not affect the link type of its closure. Therefore we may assume
$$\beta=\sigma_1^{a_1} \sigma_2^{b_1} \sigma_1^{a_2} \sigma_2^{b_1} \ldots \sigma_1^{a_m} \sigma_2^{b_m}$$
with all $a_i \geq 2$, $b_k \geq 1$. We use 
$$(a_1,\underbrace{0,\ldots,0}_{b_1-1},a_2,\underbrace{0,\ldots,0}_{b_2-1},a_3,0\ldots 0,a_m)$$
as a shortcut for $\beta$. In this notation the braids of the forbidden minors $T$, $E$, $Y$ read $(4,4)$, $(6,3)$ and $(3,0,3)$ (see again Figure~1). The surface $X$ cannot be defined by a $3$-braid. Nevertheless, it is a minor of the fibre surface of the $3$-braid $(2,0,2,0)$, as shown in Figure~3. Note that the brick diagram in the middle represents a Seifert surface in a straightforward way, but not a braid. Here and thereafter equality between brick diagrams means isotopy of surfaces or, equivalently, of their boundary links. The latter is easy to verify.
\begin{figure}[ht]
\scalebox{1.2}{\raisebox{-0pt}{$\vcenter{\hbox{\epsffile{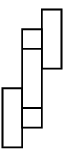}}}$}} $\quad = \quad$
\scalebox{1.2}{\raisebox{-0pt}{$\vcenter{\hbox{\epsffile{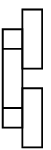}}}$}} $\quad \subset \quad$
\scalebox{1.2}{\raisebox{-0pt}{$\vcenter{\hbox{\epsffile{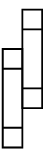}}}$}}
\caption{}
\end{figure}
We distinguish four types of braids $\beta$ depending on the number $m \in \N$.

\smallskip
If $m \geq 4$ then the braid $\beta$ (more precisely its fibre surface) contains $(2,2,2,2)$ as a minor, hence $(4,4)$, the forbidden minor $T$.

\smallskip
If $m=3$ and at least one $b_k \geq 2$ then $\beta$ contains $(2,2,2,0)$, hence $(2,0,2,0)$, in turn the forbidden minor $X$.

\smallskip
If $m=3$ and all $b_k=1$ then either $\beta$ is one of $(2,2,2)$, $(2,2,3)$ or $\beta$ contains one of $(2,3,3)$, $(2,2,4)$. The first two braids are easily identified as $E_7$ and $E_8$. The others contain the forbidden minors $(3,0,3)$ and $(4,4)$, i.e. $Y$ and $T$. 

\smallskip
If $m=2$ and both $b_k \geq 2$ then $\beta$ contains $(2,0,2,0)$, in turn $X$.

\smallskip
If $m=2$ and one $b_k \geq 2$ then either $\beta$ contains the forbidden minor $(3,0,3)$ or $\beta=(2,a,\underbrace{0 \ldots 0}_{b-1})$, which represents the same link as $(a,b+1)$. This is illustrated in Figure~4, for $a=4$, $b=5$. 
\begin{figure}[ht]
\scalebox{1.2}{\raisebox{-0pt}{$\vcenter{\hbox{\epsffile{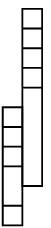}}}$}} $\quad = \quad$
\scalebox{1.2}{\raisebox{-0pt}{$\vcenter{\hbox{\epsffile{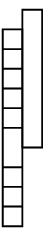}}}$}}
\caption{}
\end{figure}

\smallskip
If $m=2$ and both $b_k=1$ then either $\beta$ contains one of the forbidden minors $(6,3)$, $(4,4)$ or $\beta$ is one of $(2,n)$, $(3,3)$, $(3,4)$, $(3,5)$. These are the braids of type $D_n$, $E_6$, $E_7$ and $E_8$.

\smallskip
If $m=1$ then $\beta$ represents a $2$-stranded torus link (type $A_n$) or a connected sum of two such links.

\pagebreak

\medskip
\noindent
\textbf{$4$-braids}

\noindent
Let $L$ be a link represented by a positive $4$-braid $\beta$. We suppose that $\beta$ is irreducible and minimal, i.e. $\beta$ does not decompose into a connected sum of two non-trivial positive braids and $\beta$ has the least braid index among all positive braids representing $L$. Contrary to above, we suppose that $\beta$ has no subwords of the form $\sigma_1 \sigma_2 \sigma_1$ and $\sigma_2 \sigma_3 \sigma_2$. From this and irreducibility we deduce that $\beta$ contains the fibre surface of the braid $\sigma_2^2 \sigma_1 \sigma_2^2 \sigma_1$ as a minor. In addition, $\beta$ contains at least two non-consecutive generators $\sigma_3$. Keeping in mind that there is no subword of the form $\sigma_2 \sigma_3 \sigma_2$, we are left with one of the four minors shown in Figure~5.

\begin{figure}[ht]
\scalebox{1.2}{\raisebox{-0pt}{$\vcenter{\hbox{\epsffile{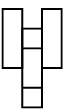}}}$}} $\quad$ 
\scalebox{1.2}{\raisebox{-0pt}{$\vcenter{\hbox{\epsffile{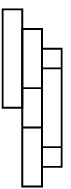}}}$}} $\quad$  
\scalebox{1.2}{\raisebox{-0pt}{$\vcenter{\hbox{\epsffile{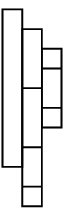}}}$}} $ = $ 
\scalebox{1.2}{\raisebox{-0pt}{$\vcenter{\hbox{\epsffile{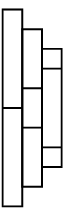}}}$}} $\quad$
\scalebox{1.2}{\raisebox{-0pt}{$\vcenter{\hbox{\epsffile{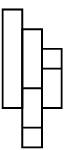}}}$}} $ = $
\scalebox{1.2}{\raisebox{-0pt}{$\vcenter{\hbox{\epsffile{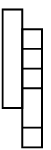}}}$}}
\caption{}
\end{figure}

The first three contain $X$ as a minor (for the third one, this is best seen after the deformation shown in the figure). The last one can be reduced to a $3$-braid and does not contain $X$ as a minor. In that case $\beta$ must have at least one more crossing of type $\sigma_2$ or $\sigma_3$ in order to be a minimal braid. A careful consideration brings up two more minors containing $X$, see Figure~6.

\begin{figure}[ht]
\scalebox{1.2}{\raisebox{-0pt}{$\vcenter{\hbox{\epsffile{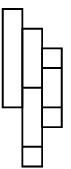}}}$}} $\quad$  
\scalebox{1.2}{\raisebox{-0pt}{$\vcenter{\hbox{\epsffile{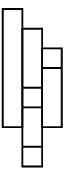}}}$}} 
\caption{}
\end{figure}

\medskip
\noindent
\textbf{$n$-braids} ($n \geq 5$)

\noindent
Let $L$ be a link represented by a positive, irreducible and minimal $n$-braid $\beta$ ($n \geq 5$).
As in the case of $4$-braids, we may suppose that $\beta$ contains a minor of the form $\sigma_2^2 \sigma_1 \sigma_2^2 \sigma_1$, as well as a symmetric version thereof on the right, $\sigma_{n-2}^2 \sigma_{n-1} \sigma_{n-2}^2 \sigma_{n-1}$. Since $\beta$ is irreducible, there exists a chain of small rectangles, successively linked, connecting a rectangle in the second column to a rectangle in the second to last column, as illustrated in Figure~7. The resulting surface manifestly contains $X$ as a minor.

\begin{figure}[ht]
\scalebox{1.2}{\raisebox{-0pt}{$\vcenter{\hbox{\epsffile{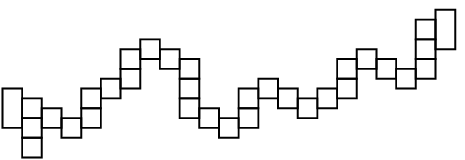}}}$}} 
\caption{}
\end{figure}

In summary, prime positive braid links of minimal braid index $n \geq 4$ have non-positive Seifert form. Non-positivity of the Seifert form is detected by four forbidden minors $T$, $E$, $X$, $Y$. At last, every prime positive braid link with positive Seifert form corresponds to a simply laced Dynkin diagram. 

It is conceivable that Theorem~1 carries over to larger classes of fibre surfaces, as long as they share certain features with positive braid surfaces. More explicitely, we may ask whether Theorem~1 holds for all fibre surfaces supporting the tight contact structure on $S^3$.

\bigskip
\noindent
Universit\"at Bern, Sidlerstrasse 5, CH-3012 Bern, Switzerland

\bigskip
\noindent
\texttt{sebastian.baader@math.unibe.ch}


\begin{thebibliography}{99}

\bibitem{AC}
     N.~A'Campo: \emph{Sur les valeurs propres de la transformation de Coxeter}, Invent. Math.~\textbf{33} (1976), no.~1, 61-67.


\bibitem{C}
     P.~R.~Cromwell: \emph{Positive braids are visually prime}, Proc. London Math. Soc.~(3) \textbf{67} (1993), no.~2, 384-424.

\bibitem{L}
     L.~Lov\'{a}sz: \emph{Graph minor theory}, Bull. Amer. Math. Soc.~\textbf{43} (2006), no.~1, 75-86.

\bibitem{N}
     T.~Nakamura: \emph{Positive alternating links are positively alternating}, J.~Knot Theory Ramifications~\textbf{9} (2000), no.~1, 107-112. 

\bibitem{Ra}
     J.~Rasmussen: \emph{Khovanov homology and the slice genus}, Invent. Math.~\textbf{182} (2010), 419-447.

\bibitem{Ru}
     L.~Rudolph: \emph{Nontrivial positive braids have positive signature}, Topology~\textbf{21} (1982), no.~3, 325-327.

\bibitem{St}
     J.~Stallings: \emph{Constructions of fibred knots and links}, Algebraic and geometric topology, Proc. Sympos. Pure Math.~\textbf{32} (1978), 55-60, Amer. Math. Soc., Providence, R.I.

\end{thebibliography}
\end{document}